\documentclass[aip,amsmath,amssymb,reprint]{revtex4-1}% 
\usepackage{times}
\usepackage{graphics}
\usepackage{bm}
\usepackage{subfigure}
\usepackage{graphicx}
\usepackage{amsthm}
\usepackage{notes2bib}
\usepackage{amssymb}
\usepackage{url}
\usepackage{enumerate}
\usepackage{epstopdf}
\usepackage{mathtools}          
\usepackage{setspace}
\usepackage[usenames,dvipsnames,svgnames,table]{xcolor}
\usepackage{color}

\newcommand{\allblack}{\color{black}{}}

\usepackage{lineno}
\begin{document}

\title{
Constructing differential equations using only a scalar time-series about continuous time chaotic dynamics
}
\author{Natsuki Tsutsumi}
\affiliation{Faculty of Commerce and Management, Hitotsubashi University, Tokyo 186-8601, Japan}
\author{Kengo Nakai}
\affiliation{Faculty of Marine Technology,  Tokyo University of Marine Science and Technology, Tokyo 135-8533, Japan}
\author{Yoshitaka Saiki}
\affiliation{Graduate School of Business Administration, Hitotsubashi University, Tokyo 186-8601, Japan}
\email{yoshi.saiki@r.hit-u.ac.jp}

\begin{abstract}
We propose a simple method of constructing a system of differential equations of chaotic behavior based on the regression only from a scalar observable time-series data.
The estimated system enables us to reconstruct invariant sets and statistical properties as well as to infer short time-series.
Our successful modeling relies on the introduction of a set of Gaussian radial basis functions to capture local structure.
The proposed method is used to construct a system of ordinary differential equations whose orbit reconstructs a time-series of a variable of the well-known Lorenz system as a simple but typical example. 
A system for a macroscopic fluid variable is also constructed.  
\end{abstract}

\maketitle

Time series data governed by some deterministic law are everywhere around us. 
However, in most cases we do not know the governing law
which enables us to predict behavior in the future.
Hence, the identification of the law, a dynamical system model which governs the observable behavior is one of the primary issues in science and technology.
In this paper we propose a simple but sophisticated method of modeling a dynamical system as differential equations only from a scalar time-series, when the dynamics is governed by a continuous time chaotic system.
The method we propose is applicable even when the observable number of variables is limited, which is often the case in the real world application. 
Our method is expected to be used to model a complex behavior in various fields such as climate science, mechanical engineering and biology, especially when a governing law is unknown. 

\vspace{-5mm}

\section{introduction}
Since the 1980s techniques of chaotic time-series analysis have been developed, and they are used to analyze experimental and observation data.
When the number of observable time series data is limited or less than the dimension of the corresponding dynamics 
as is often the case in the real world application, the delay embedding method~\cite{takens_1981,sauer1991} can be used for the  reconstruction of an attractor.
The Lyapunov spectrum, the degree of stability, is also computed from observable time series data~\cite{sano85,wolf85}.
From a time series we can obtain the existence of nonlinearity behind the background dynamics using the method of surrogate data~\cite{theiler92}.
Using the Gaussian radial basis functions time series prediction has been attempted~\cite{powell87,casdagli89}.

\allblack

We are often eager to identify a dynamical system model which generates an observable time series data. 
A dynamical model enables us to predict a behavior in the near future. Furthermore, by using the model we can recognize a rare event which has not been occurred before.
However, it is often hard to derive a governing system of complex dynamics. 
In particular, derivation of a governing system for macroscopic dynamics has been one of the fundamental unsolved problems  because of the existence of nonlinear interactions in microscopic dynamics.\\
\indent Recently several approaches have been proposed concerning modeling dynamics from given time-series data using machine learning.
These modern approaches employ high performance computations.
Neural ODE~\cite{chen2018} recently attracts much attention for using a neural network to model ODE. 
Reservoir computing is a recurrent neural network 
that is widely used for time-series inference~\cite{Pathak_2017,nakai_2018,nakai_2020,kobayashi2021}. 
The method employs relatively large number of nodes in compared to the dimension of the original dynamical system.
It does not employ differential equations. 
\allblack
{
Some studies~\cite{brunton_2017, yuan_2019}
estimate ODEs having essential structures using  simple polynomial terms behind time series.
}
Champion et al.~\cite{champion19} recently proposes a method which simultaneously learns the governing equations and the associated coordinate system by using deep neural networks and sparse identification of nonlinear dynamics. 
Their focuses are in finding the essential 
dynamical system with the fewest terms necessary to describe the dynamics behind a given data. 
They succeed in clarifying a mathematical skeleton behind the dynamics.
Model variables in their approach are not physically understandable. \\
\indent There are simple approaches which model a system of differential equations using an observable variable without using neural networks.
Baake et al.~\cite{baake1992} estimates a system of ODE by polynomial fittings. 
A higher order polynomial estimations of ODE system are attempted using a regression with Lasso regularization~\cite{wang2011,brunton16}.
These methods can be applicable only 
under the condition that 
the time derivative of each variable is approximately equal to the low order polynomials of variables.
When we have no knowledge of the choice of variables of the governing system, 
the above condition will not be satisfied.\\
\indent We propose a simple method of constructing a system of
differential equations~(ODEs) of chaotic behavior based on the regression only from
observable time-series data,
which is applicable even when the time derivative of each variable 
is not approximated by the low order polynomials of variables.
Our approach to a model construction has various advantages: (i) a model variable is physically understandable; 
(ii) a model construction requires simple steps without using a neural network; (iii) a model can be constructed even when the number of observable variables is limited; (iv) a model can be constructed even when no knowledge of the governing system is given. 
The method enables us to construct a model that is formulated explicitly using physically understandable variables.
Hence we can identify the invariant sets such as the fixed points and the periodic orbits.\\
We assume there exists an unknown system of $N$ dimensional ODEs~
\footnote{Almost all the deterministic chaotic time series of continuous time are governed by a differential equation. Such a differential equation can be autonomous differential equation, non-autonomous differential equation, partial differential equation, delay differential equation.
They can be approximated well by ordinary differential equations if they have finite dimensional attractors~\cite{temam12}.} called an original system concerning an unknown variable $\bm{x}$:
\begin{equation}
    \frac{d\bm{x}}{dt}=\bm{f}(\bm{x}).
    \label{eq:odegeneral}
\end{equation}

\indent
We can observe some or all of the components of the variable $\bm{x}$, or more generally 
\begin{equation}
    w_i=g_i(\bm{x}),~ i=1,\ldots, I. \label{eq:observable} 
\end{equation}
\begin{figure}[tb]
    \begin{center}
        \includegraphics[width=1.0\columnwidth,height=0.8\columnwidth]{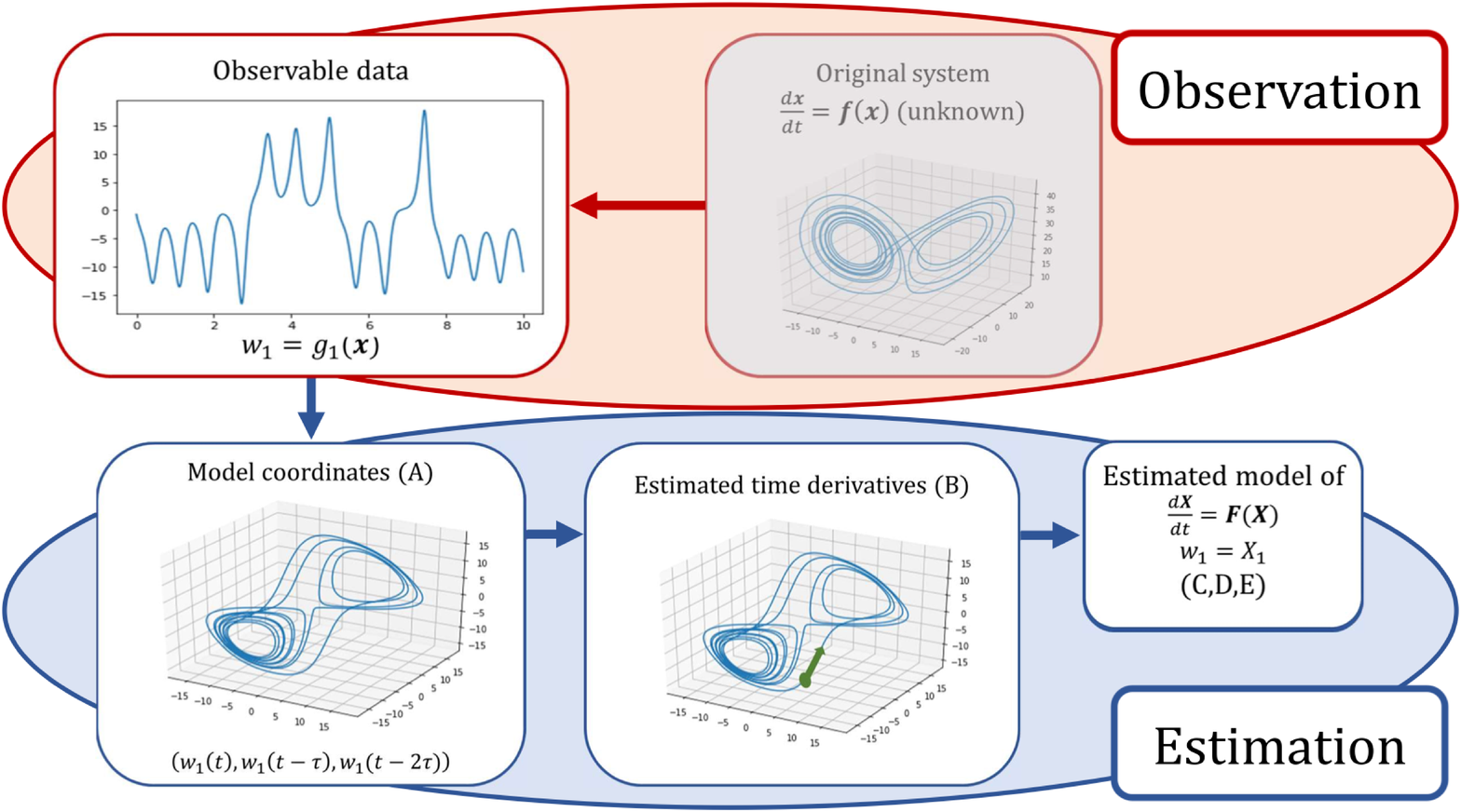}  
    \end{center}
    \caption{
    {\bf Outline of the proposed method for constructing a system of ODEs.} 
    The alphabets~(A,B,C,D,E) correspond to the sequence of steps.
    }
    \label{fig:outline}
\end{figure}
We are unable to reconstruct the original system Eq.~\eqref{eq:odegeneral} itself from given time-series data unless all the components of the variable $\bm{x}$ in Eq.~\eqref{eq:odegeneral} and their time-series are known. 
Hence, we assume that there exists a system of $D$ dimensional ODEs:
\begin{equation}
    \frac{d\bm{X}}{dt} = \bm{F}(\bm{X}),
    \label{eq:odeestimation} 
\end{equation}
where the first $I$ components of the variable $\bm{X}$ are $X_i=w_i$ $(i=1,~\ldots,I)$ and the rest of $X_i$ $(i=I+1,~\ldots,D)$ are created from the delay-coordinates of some $w_i$~\cite{takens_1981,sauer1991}.
The aim of our proposed method is to model~Eq.~\eqref{eq:odeestimation}, which is a dynamical system that can describe the behavior of the observable variables $w_i$ $(i=1,~\ldots,I)$ as components of the variable $\bm{X}$.
Throughout the paper to simplify explanations the number of observable variables $I$ is fixed as $1$ which is less than the effective dimension of the dynamics to be modeled.
In this paper, we focus on cases where the right-hand side of Eq.~\eqref{eq:odeestimation} cannot be described by low-order polynomials and present two examples.
In the first example, we construct a system of differential equations which describes the dynamics of an $x$ variable of the Lorenz system~\cite{lorenz_1963}:
    ${dx}/{dt}=10 (y-x), {dy}/{dt}=28x-y-xz, {dz}/{dt}=xy-{8z}/{3},$
when the observable time-series is limited to that of the variable $x$.
We call the Lorenz system the original system corresponding to~Eq.~\eqref{eq:odegeneral}, and a trajectory  created from it is called an actual trajectory.
At least three variables are required to model a chaotic attractor of autonomous ODEs.
Then we construct a system of ODEs of such $D=3$ dimensional variables. 
Remark that we have no knowledge of the appropriate value of the dimension $D$ and the degree of nonlinearity. 
In the second example, we construct a fluid dynamics model.
We apply the same proposed method to construct an unknown system of ODEs that governs the time development of a macroscopic fluid variable~\cite{nakai_2018}.\\
\indent 
The rest of this paper is organized as follows. 
In Section~\ref{sec:method}, we present the proposed method for deriving a system of ODEs.
In Section~\ref{sec:lorenz}, we apply the proposed method to the first example, where the only observable is the first variable of the Lorenz system. 
In Section~\ref{sec:BASIN}, we analyse the basin of attraction of the constructed model in Section~\ref{sec:lorenz}. 
In Section~\ref{sec:fluid}, we model a system of ODEs describing a macroscopic fluid flow for which the governing equations are unknown.
Concluding remarks are given in Section~ \ref{sec:conclusion}.
\allblack
\section{Method}\label{sec:method}
As shown in Fig.~\ref{fig:outline}, our aim is to construct a system of $D$ dimensional ODEs~Eq.~\eqref{eq:odeestimation} based only on some observable deterministic time-series of discrete time.
The steps of the proposed method are outlined below:
\begin{itemize}
    \item[{\textsf{A.}}] Choose the delay-coordinates: delay-time and dimension 
    \item[{\textsf{B.}}] Estimate the time derivative at sample points using the Taylor expansion
    \item[\textsf{C.}] Choose the basis function\footnote{ We do not assume some particular vector space in advance as is usual in the literature for machine learning.} used in Step \textsf{D} 
    \item[\textsf{D.}] 
    Perform ridge regression at sample points
    \allblack
    \item[\textsf{E.}] { Evaluate the model quality according to the reproducibility of
    a 
    \allblack
    delay structure.}
\end{itemize}

\subsection{Case of limited variables}
\label{sec:delay}
For practical purposes, as the number of observable variables is small, compared with the dimension of the background dynamical system. 
We can use the delay-coordinates $(w_1(t),w_1(t-\tau), w_1(t-2\tau),\ldots, w_1(t-(D-1)\tau))$ of a limited observable  to obtain a system of ODEs.
See Sec.~\ref{sec:lorenz} for the choice of $\tau$.

\subsection{Estimating the time derivative through smoothing the trajectory}
\label{sec:timederivative}
We can use the sixth-order Taylor expansion to estimate the time derivative at each sample point $\bm{X}(\tilde{t})$ based only on discrete time points of a trajectory $\bm{X}(t)$ $(t=\ldots, \tilde{t}-\Delta t,\tilde{t},\tilde{t}+\Delta t,\ldots)$:
\begin{align*}
    \frac{d\bm{X}(\Tilde{t})}{dt} \approx 
    \{ \bm{X}(\Tilde{t}+3\Delta t) - 9\bm{X}(\Tilde{t}+2\Delta t) +45 \bm{X}(\Tilde{t}+\Delta t) \\
    - 45\bm{X}(\Tilde{t}-\Delta t) + 9\bm{X}(\Tilde{t}-2\Delta t) - \bm{X}(\Tilde{t}-3\Delta t)\}/(60 \Delta t),
\end{align*}
where $\Delta t$ is the time step of sample points along a training trajectory.
Lower- and higher-order Taylor expansions will also work well.
{%\allblue 
When the observable data includes noise, we can use $l\Delta t$ for some positive integer $l$ instead of $\Delta t$
to estimate the time derivative, which enables to avoid high frequency oscillations.
}

\subsection{Choice of a basis function}
\label{sec:basefunctions}
When we model Eq.~\eqref{eq:odeestimation} by linear regression using limited computational resources, 
basis functions should be chosen appropriately.
After we explain polynomial basis functions
used in Wang et al.~\cite{wang2011}, 
we introduce Gaussian radial basis functions. 
We combine them with polynomial basis functions for the basis of our regression. 

\subsubsection{Polynomial basis function}
For ODE estimation through regression, one possible choice is the polynomial basis function: 
\begin{equation}
    F_k(\bm{X}) \approx \tilde\beta_0^k + \sum_{d=1,\cdots,D} \tilde\beta_d^k X_d +\sum_{1 \leq i \leq j \leq D}
    \tilde\beta_{e(i, j)}^k X_{i}X_{j}+\cdots,
    \label{eq:polynomial}
\end{equation} 
where 
$e(i, j) = iD+j-\frac{1}{2}(i-1)i$ and  $\tilde{\boldsymbol{\beta}}{^k}$ is a vector  whose components are $\tilde\beta^k_l$ $(l=0,1,2,\ldots)$, and  $F_k(\bm{X})$ is the $k$th component of $\bm{F}(\bm{X})$. 

Wang et al.~\cite{wang2011} employ the polynomial basis function to estimate dynamical systems for the H\'enon map, R\"ossler system and Lorenz system based on time-series data. They report that bifurcation diagrams are reconstructed well. They used Lasso regularization for the regression of the  coefficients.
In this paper we assume no knowledge of variables of a governing equation. Furthermore, 
the function $\bm{F}$ in Eq.~\eqref{eq:odeestimation}
 generally does not have the form of low-order polynomials.
We confirmed that the 
method is not applicable to our problem of 
ODE construction
(see Appendix~\ref{appendix:polynomial}).
 Although it tends to be successful at short-term time-series inference, 
 it tends to fail at reconstructing statistical properties.

\subsubsection{Gaussian radial basis function}
\begin{figure}[t]
    \vspace*{2mm}
    \begin{center}
        \includegraphics[width=0.98\columnwidth,height=0.42\columnwidth]{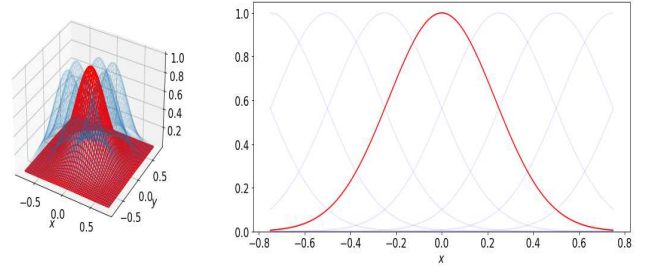}
    \end{center}
    \caption{{\bf Gaussian radial basis function $\phi_j$ in Eq.~\eqref{eq:gaussian}.}  
    On the left 2D function centered on $(0,0)$~(red) together with functions centered on $\{c_j\}=(0.25n_x,0.25n_y)~(n_x,n_y=-1,0,1)$ (blue) is plotted. 
    On the right 1D function centered  on $0$~(red) is plotted together with functions centered on $\{ c_j \}=0.25n_x~(n_x=-4,\ldots,4)$ (blue). 
    }
    \label{fig:gaussianradialfunc}
\end{figure}

To overcome the difficulty of modeling local structures in complex dynamics, 
we can approximate $\bm{F}(\bm{X})$ by using 
the localized Gaussian radial basis function $\phi_j$:
\begin{equation}
F_k(\bm{X}) \approx \tilde\beta_0^k + \sum_{d=1,\cdots,D} \tilde\beta_d^k X_d + \sum_{j=1,\cdots,J} \tilde\beta_{D+j}^k~\phi_j(\bm{X}),\label{eq:gaussianpoly}
\end{equation}
where $\tilde{\boldsymbol{\beta}}{^k}$ is a set of estimated parameters and  
\begin{equation}
\phi_j(\bm{X}) = \exp\left(\frac{-||\bm{X}-c_j||^2}{\sigma^2}\right)\label{eq:gaussian},
\end{equation}
where $\|\cdot\|$ denotes the $l^2$ norm.   %footnote? to tsutusmi
Here $c_j$ is the coordinate of the $j$~th center point ($j=1,\ldots,J$), and $\sigma^2$ is the parameter 
that determines the deviation of $\phi_j$.
The parameter $c_j$ is distributed as lattice points with the grid size $\delta_{grid}$.
In this paper, $\sigma^2$ is approximately equal to $1.7372 \delta_{grid}^2$
\hspace{-2mm}\footnote{ 
    The choice of $\sigma^2$ is  obtained from  
    $$
        \sigma^2 := \frac{((m-1) \delta_{grid})^2}{- \log_{e} p}, 
    $$
    where $m$ is the degree of the corresponding B-spline basis function and $p~(>0)$ is a small value. We consider $c_j$ only if there exists a data point in the $(m-1)\delta_{grid}$-neighborhood.
    We set $m=3$, and $p=0.1$, then 
    $(m-1)^2 /(- \log_{e} p) \approx 1.7372$.}
    \hspace{-1mm}.
Figure~\ref{fig:gaussianradialfunc} shows the shape of $\phi_j$ in Eq.~\eqref{eq:gaussian} for $D=1$ or $2$ with $\delta_{grid}=0.25$. 
See Kawano and Konishi~\cite{kawano2007} for more details.

%\allblue
\subsection{Ridge regression}
\label{sec:regression}
\allblack
To estimate a model of Eq.~{\eqref{eq:odeestimation}}, we can apply least squares regression with an appropriate basis function using  ridge regularization. 
Here we briefly review the regression method.
We can estimate the coefficients $\boldsymbol{\beta}$ of a linear model $y_i\approx  \bm{a}_i^\mathrm{T} \cdot \boldsymbol{\beta}$, where $y_i$ is a normalized target variable,  $\bm{a}_i$ is a set of normalized explanatory variables and ${}^\mathrm{T}$ represents a transpose of the matrix. 
For this model, the ridge estimator is realized by minimizing the following function $L(\bm{b})$:
\begin{align}
    L(\bm{b}) = \frac{1}{2n} ||\bm{y} - A \bm{b} ||^2 + \frac{\lambda}{2} ||\bm{b}||^2,\label{eq:regularization}
\end{align}
where $n$ is the size of data, $\bm{y} = [y_1, y_2, \ldots, y_n]^\mathrm{T},$ $A=[\bm{a}_1,\bm{a}_2, \ldots ,\bm{a}_n]^\mathrm{T}$, and $\lambda$ is a parameter determining the strength of regularization.
The ridge estimator can be written as follows: 
\begin{equation*}
    \boldsymbol{\beta} = (A^\mathrm{T} A + n \lambda I )^{-1} A^\mathrm{T}  \bm{y}, 
\end{equation*}
where $I$ is the identity matrix. 
When we estimate the system of ODEs~Eq.~\eqref{eq:odeestimation}, $\bm{y}$ represents the time derivative, and $A$ corresponds to  $\{ x_{(1,i)}, \ldots x_{(D, i)}, \phi_1(\bm{X}_i), \ldots, \phi_J(\bm{X}_i) \}_{i=1, \ldots, n}$
(see \eqref{eq:gaussianpoly} in Sec.~\ref{sec:basefunctions}).

\subsection{Evaluation of the model quality}
\label{sec:delaystructure}
For the model to be appropriate we can confirm that 
the delay structure $(w_1(t),w_1(t-\tau), w_1(t-2\tau),\ldots, w_1(t-(D-1)\tau))$ is %are 
reconstructed in the model trajectory.
See Figs.~\ref{fig:delay} and \ref{fig:delaydistribution} in Section~\ref{sec:lorenz} and Fig.~\ref{fig:fluidshort} in Section~\ref{sec:fluid}.
The hyper-parameters such as regularization parameter and grid size are chosen appropriately based on the degree of the  reconstruction.

\allblack

\section{Construction of data-driven ODEs: Lorenz dynamics}\label{sec:lorenz}
As the first example we construct a system of ODEs that describes the behavior of the Lorenz system.
We assume the number of observable variables in \eqref{eq:observable} is $I=1$, and the variable is $w_1=x,$ where $x$ is a variable of the original Lorenz system.  
We set the model coordinate $\bm{X}$ in Eq.~\eqref{eq:odeestimation} as  $\bm{X}(t)=(w_1(t),w_1(t- \tau),w_1(t-2 \tau))$, where $ \tau=0.13$ ($D=3$ in Eq.~\eqref{eq:odeestimation}).
\footnote{
    The auto-correlation is referred to for choosing the delay-coordinate.In this case the correlation coefficient of $w_1(t)$ and $w_1(t- \tau)$ is 0.79, and that of $w_1(t)$ and $w_1(t-2 \tau)$ is 0.46.
}
The time length $T$ of the training time-series is $5000$ (time step $\Delta t=0.005$), and $2\%$ are chosen as the sample points for regression. 
The number of center points $c_j (j=1,\ldots,J)$  of the Gaussian radial basis function is $1806$ with a grid size of about $2$ corresponding to the product of the standardization coefficient and normalized grid size~($\approx7.9261\cdot 0.25$). 
The regularization parameter $\lambda=10^{-7}.$ 
\\
\indent We evaluate a model from various aspects.
After confirming that the constructed model has a trajectory approximate to that created from the actual model, we investigate invariant sets  such as the chaotic attractor and fixed points. 
We also confirm the reconstruction of the delay structure among components among a variable $\bm{X}$.\\

\medskip

\begin{figure}[tb]
    \vspace*{2mm}
    \begin{center}
        \includegraphics[width=0.98\columnwidth,height=0.5\columnwidth]{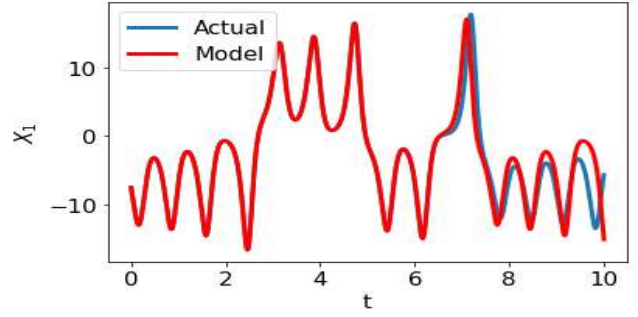}
    \end{center}
    \caption{{\bf A short-term trajectory of $X_1$ of the model in compared to the corresponding trajectory of $x$ of the original Lorenz system.}
    These two trajectories behave similarly but deviate after some time because of the chaotic property.
    }
    \label{fig:shortorbits}
\end{figure}

\indent {\bf Short-term orbit.}
We found that a time-series inference of $x$ can successfully be applied for a short time 
under many initial conditions.
Figure~\ref{fig:shortorbits} shows an example trajectory. 
The long-term time-series inference inevitably fails because of the chaotic property of the Lorenz system. 
\begin{figure}[t]
    \vspace*{2mm}
    \begin{center}  
        \includegraphics[width=0.98\columnwidth,height=0.5\columnwidth]{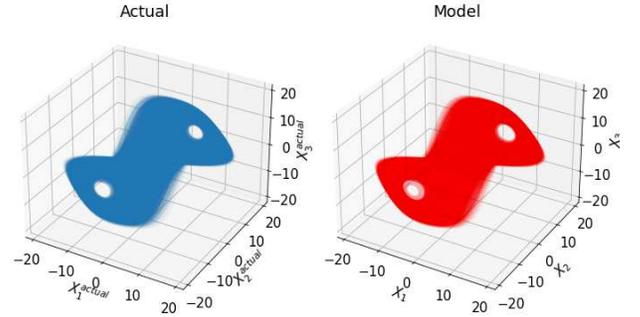}
    \end{center}
    \caption{{\bf Projections of long-term  trajectories.} A long-term trajectory of the model approximates the attractor of the original Lorenz system projected onto the $(X_1,X_2,X_3)$ space.
    We use time length $T=10,000$. 
    }
    \label{fig:attractors}  
\end{figure}

\medskip

\begin{figure}[t]
    \vspace*{2mm}
    \begin{center}
        \includegraphics[width=0.95\columnwidth,height=0.6\columnwidth]{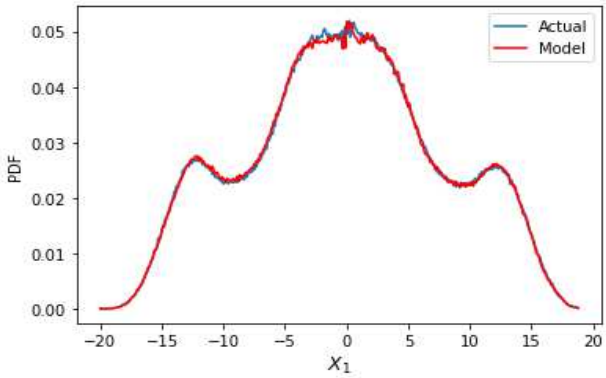}  
    \end{center}
    \caption{{\bf Density distribution of $X_1$ of the model in compared to that of the original Lorenz system.} The difference in area between the two distributions is {0.00418.}
    Each distribution is computed from a single trajectory with time length $T=10,000$. 
    }
    \label{fig:density}
\end{figure}

{\bf Attractor and invariant density.}~We find it difficult to represent long-term behavior accurately from a model constructed  using only polynomial basis functions alone for the regression~(see Appendix~\ref{appendix:polynomial}).
Here, we show that a model constructed using the proposed method can reconstruct a statistical quantity.
Figure~\ref{fig:attractors} shows long-term trajectories.  
Figure~\ref{fig:density} shows an agreement of the density distributions computed from a model trajectory $\{X_1(t)\}$ and a trajectory $\{x(t)\}$ of the original Lorenz system. 
\begin{figure}[tb]
    \vspace*{2mm}
    \begin{center}
        \includegraphics[width=0.95\columnwidth,height=0.6\columnwidth]{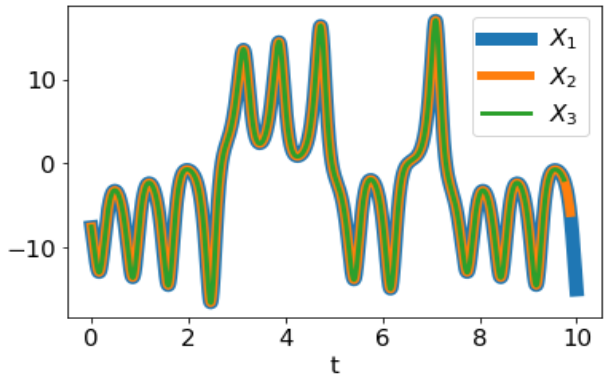}
    \end{center}
    \caption{{\bf Reconstruction of time delay among components $X_1$, $X_2$ and $X_3$ of a variable $\bm{X}$.} 
    The time-series of $X_1(t), X_2(t+\tau)$ and $X_3(t+2\tau)$ ($0 \leq t \leq 10$) of the same trajectory as that in  Fig.~\ref{fig:shortorbits} are shown to be close to each other,  which confirms that the delay structure among the three variables is reconstructed by the model. Even though the model orbit deviates from the orbit of the original system around $t=7$ in Fig.~\ref{fig:shortorbits}, the delay structure remains satisfactory, which implies the validity of the model.
    }
    \label{fig:delay}
\end{figure}

The Lyapunov exponents are used to evaluate the degree of instability.
We find that the positive and zero exponents of the model are  $\tilde{\lambda}_1=0.898$ and $\tilde{\lambda}_2=0.000$, which  agree well with those of the original Lorenz system~($\lambda_1=0.906$ and $\lambda_2=0.000$). The third negative exponent depends on the choice of coordinates; thus, the disagreement between that of the model ($\tilde{\lambda}_3=-5.730$) and that of the original Lorenz system ($\lambda_3=-14.572$)) is reasonable. 

\medskip

\begin{figure}[tb]
    \vspace*{2mm}
    \begin{center}
        \includegraphics[width=0.95\columnwidth,height=0.43\columnwidth]{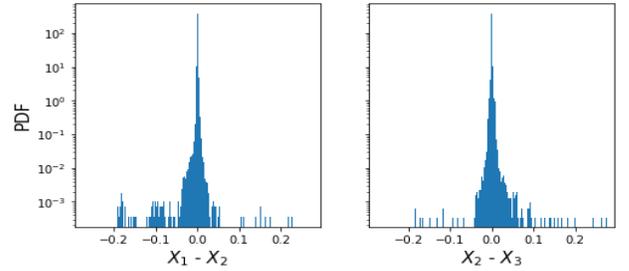}
    \end{center}
    \caption{{\bf Density distributions of $X_1(t)-X_2(t+ \tau)$~(left), and of $X_2(t)-X_3(t+ \tau)$~(right).}
    The vertical axis is on a logarithmic scale. 
    Each distribution is computed from a trajectory with time length $T=10,000$.
    The standard deviations of the distributions are $0.0015$~(left) and $0.0015$~(right);  each is localized around zero  with the standard deviation $\sigma_{X_1}$ of $X_1$ ($\sigma_{X_1} =7.9261$).
    This confirms that the delay structure among components of a variable is reconstructed by the model.
    }
    \label{fig:delaydistribution}  
\end{figure}

{\bf Delay structure.}
As we employ the delay-coordinate  with $D=3$ and $\tau=0.13$ for the variable $\bm{X}$ of a model, 
the relation 
    $X_1(t)\approx X_2(t+ \tau) \approx X_3(t + 2 \tau)$
should hold for a model to be appropriate.
Figure~\ref{fig:delay} shows time-series $X_1(t)$, $X_2(t + \tau)$ and $X_3(t + 2 \tau)$ that satisfy the relation. 
Figure~\ref{fig:delaydistribution} confirms that the density distributions of $X_1(t)-X_2(t+ \tau)$ and  $X_2(t)-X_3(t+ \tau)$ are each localized around zero, which suggests the successful reconstruction of the delay structure.
    Note that the delay structure is not reconstructed well when a ridge  regularization parameter for the regression is not chosen appropriately (see Appendix~\ref{appendix:delaystructure}).
\allblack

\medskip

{\bf Fixed points.}
Fixed points are the most fundamental invariant sets.
The original Lorenz system has three fixed points,  whose coordinates are $(\pm8.48528,$ $ \pm8.48528, \pm8.48528)$ and $(0,0,0)$ in the delay-coordinate.
The model has the fixed points $L_{\text{model}}$, $R_{\text{model}}$ and $O_{\text{model}}$, as given in Table~\ref{tab:fixed-points}. 
The model has two additional fixed points, that we call the ghost fixed points $GL_{\text{model}}$ and $GR_{\text{model}}$ outside the attractor of the model, and their coordinates strongly depend on the value of the regularization parameter of ridge regression. 
The role of the ghost fixed points in relation to the basin of attraction of the data-driven model 
are discussed in the following section.
\allblack
\begin{table}[t]
    \begin{center}
    \small
	\begin{tabular}{|l|r|r|r|r|r|r}
	\hline  	
    & $L_{\text{model}}$& $R_{\text{model}}$&
           	     $O_{\text{model}}$& $GL_{\text{model}}$& $GR_{\text{model}}$ \\ \hline
    $x_1^*$ & $-8.4788$ & $8.4786$ & $0.0009$ & $-1.2879$ & $1.3244$  \\ \hline
    $x_2^*$ & $-8.4689$ & $8.4724$ & $0.0008$ & $-1.2498$ & $1.2864$  \\ \hline
    $x_3^*$ & $-8.4802$ & $8.4804$ & $0.0010$ & $-1.3695$ & $1.4097$ \\ \hline
    $\Lambda_1^*$ & $ 0.10+10.17i$ & $ 0.10+10.17i$ & $ 10.56$ & $ 11.19$ & $ 11.34$  \\ \hline
    $\Lambda_2^*$ & $ 0.10-10.17i$ & $ 0.10-10.17i$ & $  2.83$ & $- 5.12$ & $- 4.78$  \\ \hline
    $\Lambda_3^*$ & $-10.60$ & $-11.22$ & $-15.76$ & $-10.12$ & $-10.62$ \\ \hline
     \end{tabular}
    \normalsize
    \caption{{\bf Coordinates $(x^*_1,x^*_2,x^*_3)$ and eigenvalues $(\Lambda^*_1,\Lambda^*_2,\Lambda^*_3)$ of each fixed point of the model. 
 		} 
 		The coordinates for the fixed points of the original system projected onto the delay-coordinate $(x(t), x(t- \tau), x(t-2 \tau))$
 		are $(\pm8.4853,\pm8.4853,\pm8.4853)$, $(0,0,0)$.
 		The three fixed points $L_{\text{model}}$, $R_{\text{model}}$, and $O_{\text{model}}$ are 
 		two-dimensionally unstable and one-dimensionally stable, 
 		whereas
 		the two ghost fixed points outside the attractor 
 		$GL_{\text{model}}$ and $GR_{\text{model}}$ 
 		are one-dimensionally unstable and two-dimensionally stable.
    }
    \label{tab:fixed-points}
    \end{center}
\end{table}

\section{Basin of the model attractor}\label{sec:BASIN}
\begin{figure}[t]
    \vspace*{2mm}
    \begin{center}
        \includegraphics[width=0.98\columnwidth,height=0.55\columnwidth]{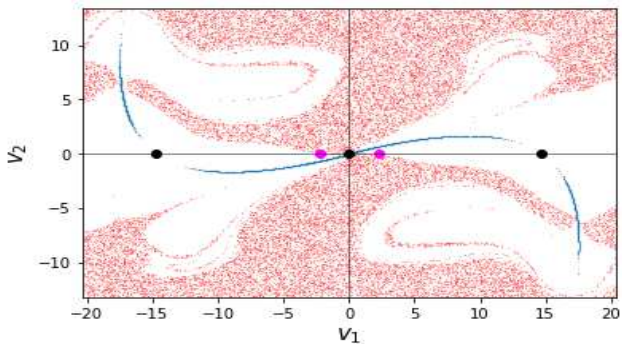}  
    \end{center}
    \caption{{\bf Basin of attraction of the data-driven model projected onto $(v_1,v_2)$ space spanned by $(1,1,1)$ and $(1,-1,0)$.} A trajectory starting from a point in a red region goes away within time $5$ and a trajectory starting from a point in a white region does not. The three fixed points in black are embedded in the attractor of the data-driven model and the two fixed points in magenta are outside the attractor and on the basin boundary.
    The set of points in blue approximates the attractor.
    }
    \label{fig:basin-of-attoractor}
\end{figure}

\begin{figure}[t]
    \vspace*{2mm}
    \begin{center}
        \includegraphics[width=0.98\columnwidth,height=0.55\columnwidth]{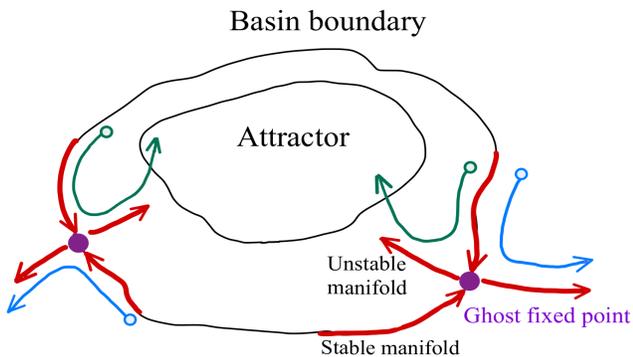}
    \end{center}
    \caption{
    {\bf Conceptual image of the model dynamics concerning the basin of attraction.} 
    The basin boundary is considered to be formed by two dimensional stable manifold of the ghost fixed points.
    }
    \label{fig:conceptual}
\end{figure}

\indent We compute the basin of attraction of the constructed data-driven model of the Lorenz system in order to capture the validity and the limitation of the model.
A set of points that are attracted to the attractor of a model is called the basin.  
The basin of attraction of the original Lorenz system is known to be $\mathbb{R}^3$. 
The white region in Fig.~\ref{fig:basin-of-attoractor} shows the basin of attraction together with the model attractor and the five fixed points of the data-driven model. 
Although three of the fixed points embedded in the attractor correspond to those of the original Lorenz system, the other two do not exist in the original system.
We find that the two fixed points $GL_\text{model}$ and $GR_\text{model}$ with two stable directions and one unstable direction in Table~\ref{tab:fixed-points} 
are on the basin boundary formed by the  
stable manifold of each of the two fixed points (See Fig.~\ref{fig:conceptual}). 
The computation of trajectories of the original system under $10^{10}$ initial conditions  confirm that, in $(X_1,X_2,X_3)$-coordinates,
at least one of the trajectories visit 
approximately $99.7\%$ of the $80^3$ boxes in $[-20,20]^3$.
We confirm that a model for $D=4$ with $ \tau=0.09$ has two ghost fixed points whose dimension of the stable manifold is three, and the same discussion will work. 

\medskip
\section{Construction of data-driven ODEs: Macroscopic variable of fluid flow}
\label{sec:fluid}
\begin{figure}[ht]
    \begin{center}
        \includegraphics[width=0.95\columnwidth,height=0.7\columnwidth]{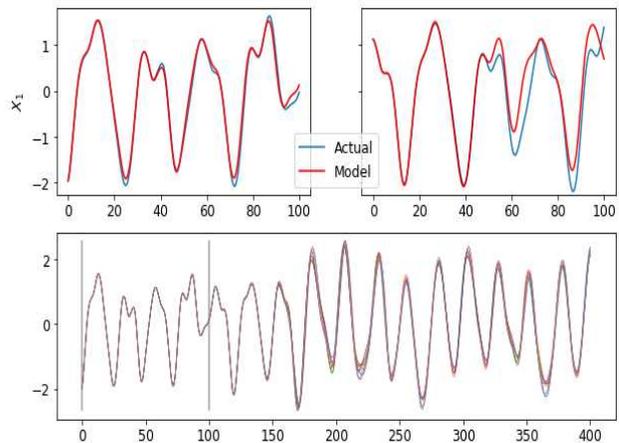}
    \end{center}
    \caption{{\bf Short-term trajectories of $X_1=E(3,t)$ of the model and reconstruction of 
    %\allblue
    the 
    \allblack
    delay structure for fluid dynamics. } 
    Short-term trajectories of $X_1=E(3,t)$ of the model under two different initial conditions and the corresponding actual trajectory are shown in the top panels. 
    The time-series of $X_i(t+(i-1)\tau)$ $(0\le t \le 400)$  is shown for $i=1,\ldots,8$ in the bottom panel. 
    The time period $0\le t \le 100$ corresponds to that of the top left panel.
    The model trajectory is shown to approximate the corresponding actual trajectory for a certain amount of time and deviate due to the chaotic property, but continues to satisfy {the} time delay structure among variables.
    }
    \label{fig:fluidshort}
\end{figure}
For the second example, we construct a system of ODEs that governs macroscopic behavior of a fluid flow, which is known to be a difficult task.
Macroscopic data of a fluid flow is obtained from a direct numerical simulation of the Navier--Stokes system.
We assume the number of observable variables in Eq.~\eqref{eq:observable} is $I=1$. The observable variable is $w_1(t)=E(3,t)$ corresponding to the energy at a certain wavenumber range which is normalized in our analysis.
See Nakai and Saiki~\cite{nakai_2018}.
We set the model variables $\bm{X}$ in Eq.~\eqref{eq:odeestimation} as  $\bm{X}(t)=(w_1(t),w_1(t- \tau),\ldots, w_1(t-7 \tau))$, where $ \tau=1.5$ ($D=8$ in Eq.~\eqref{eq:odeestimation}). 
Note that we do not have a knowledge of the Navier--Stokes equation during the construction procedure. 
The time length $T$ of the training time-series is $10^5$ (time step $\Delta t=0.05$), $2\%$ of which are chosen as the sample points for regression. 
The number of center points $c_j(j=1,\ldots,J)$ of the Gaussian radial basis function is  $922,598$ with the grid size of $0.25.$
The number of grids is restricted by our limited computational resources.
The regularization parameter $\lambda=10^{-6}.$ 
\medskip

\indent{\bf Short-term trajectories and 
%\allblue
the 
\allblack
delay structure.}
We find that time-series inference is successful for some time by integrating the obtained model.
Two panels in Fig.~\ref{fig:fluidshort} (top) confirm that a single model could infer the time-series of $E(3,t)$ under different initial conditions.
As we employ the delay-coordinate as a model variable, we expect 
$X_1(t)\approx X_2(t+ \tau)\approx \cdots \approx X_8(t + 7 \tau).$
The relation is confirmed for the short time-series of $X_1(t), X_2(t+ \tau), \cdots, X_8(t + 7\tau)$ as shown in Fig.~\ref{fig:fluidshort} (bottom).

\section{Concluding remarks}\label{sec:conclusion}
We propose a simple method of constructing a system of differential equations by using observable time-series data.
A model constructed using our method enables to reconstruct a long term dynamics, and the basin structures of the model are studied.
We exemplify that the method can be applicable for modeling  a macroscopic fluid flow as well as the chaotic Lorenz system even when the number of observable variable is limited to one. 
As a future work, it is interesting to confirm the applicability of our method to various types of chaotic dynamics including an infinite dimensional one.
\allblack

%%%%%%%%%%%%%%%%%%%%%%%%%%%%%%%%%%%%%%%%%%%%%%%%%%%%%%%%%%%%%%%%%%%%%%%%
\appendix
\section{Failure of constructing ODEs by polynomial regression}
\label{appendix:polynomial}

\begin{figure}[ht]
    \vspace*{2mm}
    \begin{center}
        \includegraphics[width=0.98\columnwidth,height=0.5\columnwidth]{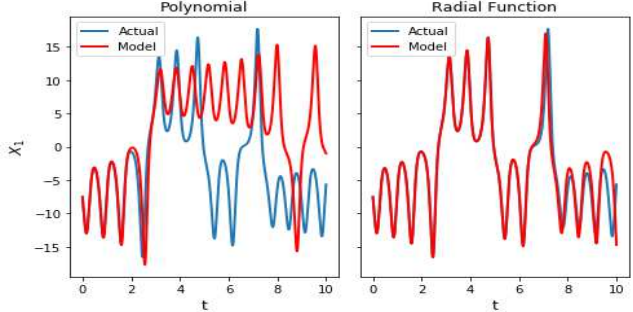}
    \end{center}
    \caption{{\bf Short-term trajectories of two models.} 
    The left and right panels show short-term  trajectories of a model by using only polynomial terms and of the main model, each of which is plotted together with the corresponding trajectory of the original system.}
    \label{fig:short-orbit_poly}  
\end{figure}

\begin{figure}[th]
    \vspace*{2mm}
    \begin{center}
        \includegraphics[width=0.98\columnwidth,height=0.5\columnwidth]{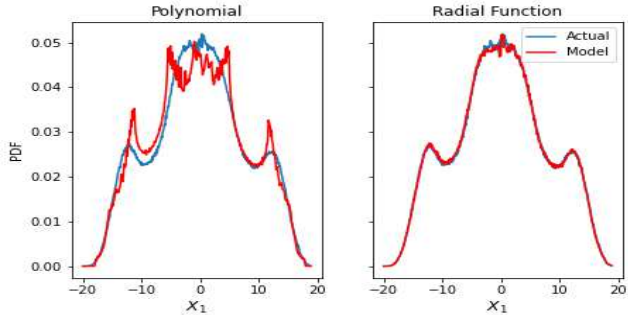}
    \end{center}
    \caption{{\bf Density distributions of $X_1$ for two models together with that of the original system: using only polynomial terms~(left) and the main model~(right).} 
    The errors in area from the distribution of the original system are 0.0235~(left) and 0.0042~(right).
    }
    \label{fig:density_poly}
\end{figure}
For the estimation of ODE in Eq.~\eqref{eq:odeestimation} through regression, one possible choice is the regression using only polynomial basis functions Eq.~\eqref{eq:polynomial}. 

Here, we present the difficulty of constructing ODEs by regression using only polynomial terms. More precisely, inferring a time-series in the  short-term is relatively easy, but the density distribution is difficult to reconstruct because of the difficulty with modeling long-term behavior. 

\begin{figure}[t]
    \vspace*{2mm}
    \begin{center}
        \includegraphics[width=0.98\columnwidth,height=0.47\columnwidth]{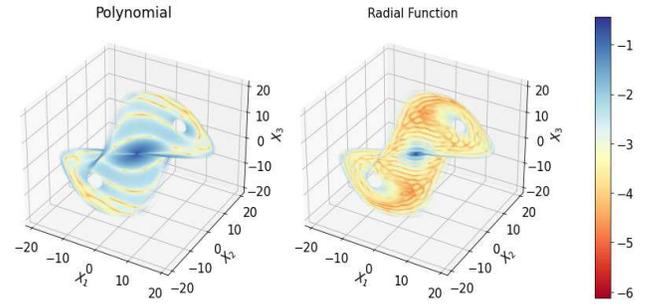}
    \end{center}
    \caption{{\bf Relative regression errors in logarithmic scale of the two models:using only polynomial terms and of the main model.} 
    The average relative errors are 0.0148~(left) and 0.0009~(right).
    }
    \label{fig:error-plot_poly}  
\end{figure}

Assume that the variable $\bm{X}(t)$ follows Eq.~\eqref{eq:odeestimation} and $\bm{F}$ is sufficiently smooth. % for some proper value of $k>0$.
We can estimate $\bm{F}$ by regression using only  polynomial basis functions Eq.~\eqref{eq:polynomial}. 
In fact, previous researchers have reported that simple systems such as the Lorenz system and R\"ossler system~\cite{baake1992,wang2011} can be successfully estimated from time-series when $\bm{F}$ is written using quadratic terms. 
However, our aim is to construct ODEs when $\bm{F}$ is not necessarily written by low order polynomials.
We apply 8th order polynomial regression to estimate a model  with the parameter 
$\lambda=10^{-6.2}$, where $\lambda$ is a coefficient of ridge regularization. 
Figures~\ref{fig:short-orbit_poly} and \ref{fig:density_poly} confirm that the model can reproduce a short orbit but cannot approximate a density distribution created from a long orbit. 
Figure~\ref{fig:error-plot_poly} shows that the relative regression error is larger for various $\bm{X}$ than for the main model. 
This implies that localized structures cannot be captured well polynomial terms up to the eight order at least.  

%\allblue 
In theory the use of higher order polynomial terms will improve the model quality. However, in practice it is difficult to construct a model because of the computational resources and the difficulty in the choice of regularization parameter in ridge regression.

%%%%%%%%%%%%%%%%%%%%%%%%%%%%%%%%%%%%%%%%%%%%%%%%%%%%%%%%%%%%%%%%%%%%%%%%%
\section{Choice of a regularization parameter from the reproducibility of the delay structure}\label{appendix:delaystructure}
\begin{figure}[ht]
    \vspace*{2mm}
    \begin{center}
        \includegraphics[width=0.98\columnwidth,height=0.55\columnwidth]{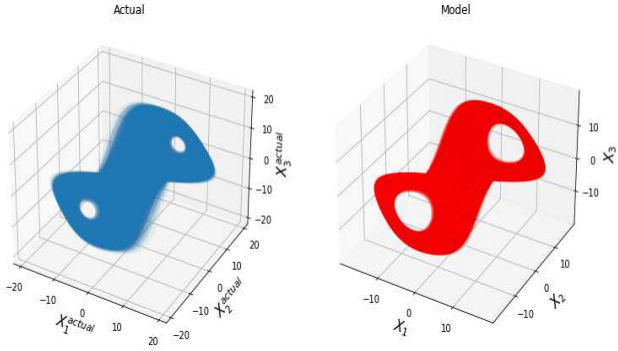}
    \end{center}
    \caption{{\bf Projection of a long-term  trajectory to $(X_1,X_2,X_3)$.} Approximation of an attractor of an constructed model with the  regularization parameter  $\lambda = 10^{-3.9}$~(right) and that of a trajectory of the original Lorenz system~(left) projected onto the delay coordinate.}
    \label{fig:attractor5.5}
\end{figure}

\begin{figure}[t]
    \vspace*{2mm}
    \begin{center}
        \includegraphics[width=0.98\columnwidth,height=0.38\columnwidth]{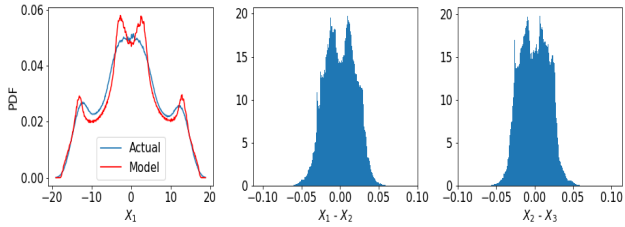}
    \end{center}
    \caption{ {\bf Density distributions for a model with the regularization parameter 
    { $\lambda = 10^{-3.9}$}.} 
    The discrepancy in density distributions of $X_1$ along the trajectories of the model and the original system~(left), and the distribution of $X_1(t)-X_2(t+\tau)$ ($\sigma=0.0196$)~(center), and that of $X_2(t)-X_3(t+ \tau)$ ($\sigma=0.0186$)~(right) of the model. 
    The discrepancy shown in the left panel is consistent with the fact that the delay relations are not satisfied well in comparison with Fig.~\ref{fig:delaydistribution} for a model with the regularization parameter $\lambda = 10^{-3.9}$.
    }
    \label{fig:delaydistribution5.5}
\end{figure}

We choose a regularization parameter in Eq.~\eqref{eq:regularization} based on the 
reproducibility of the delay structure.
In constructing a ODE system we introduce delay-coordinate variables of an observable variable when the number of observable variables is not large relative to the dimension of the background dynamical system. 
We expect that the constructed model could recover delay relations among variables.
For the first example of modeling the Lorenz dynamics, the relation $X_1(t)\approx X_2(t+\tau)\approx X_3(t + 2\tau)$ should hold for the model to be appropriate.
This relation is confirmed by the short time-series of $X_1(t)$, $X_2(t+\tau)$, $X_3(t+2\tau)$ shown in Fig.~\ref{fig:delay} and by the distributions of $X_1(t)-X_2(t + \tau)$, $X_2(t)- X_3(t +  \tau)$ having sharp peaks around zeros as shown in Fig.~\ref{fig:delaydistribution}. 
Figures \ref{fig:attractor5.5} and \ref{fig:delaydistribution5.5} show the attractor and density distributions of $X_1(t)$, $X_1(t)-X_2(t+ \tau)$ and $X_2(t)-X_3(t+\tau)$ of the model when $\lambda=10^{-3.9}$. 
These figures suggest that the model quality is worse with $\lambda=10^{-3.9}$ than with $\lambda=10^{-7},$ although the degrees of the model reconstruction are similar.
\allblack
\newpage
\section*{\bf Acknowledgements.} 
    YS was supported by JSPS KAKENHI (19KK0067, 21K18584). 
    KN was supported by JSPS KAKENHI (22K17965) and the Project of President Discretionary Budget of TUMST. 
    Part of the computation was supported by JHPCN (jh210027, jh220007), HPCI (hp210072), and the Collaborative Research Program for Young$\cdot$Women Scientists of ACCMS and IIMC, Kyoto University.

\end{document}